\documentstyle{amsppt}
\magnification1200
\pageheight{8in}
\NoBlackBoxes
\def\d{{\underline{\delta}}}
\def\L{\fracwithdelims()}

\topmatter
\title Zeros of Fekete polynomials \endtitle
\author Brian Conrey, Andrew Granville, Bjorn Poonen and K. Soundararajan 
\endauthor
\address (Conrey) American Institute of Mathematics. 360 Portage Ave, 
Palo Alto, California 94306, USA\endaddress
\email conrey\@aimath.org \endemail
\address (Conrey) Department of Mathematics, Oklahoma State University,
Stillwater, Oklahoma 74078, USA\endaddress
\email conrey\@hardy.math.okstate.edu \endemail
\address (Granville) Department of Mathematics, The University of Georgia,
Athens, Georgia 30602-7403, USA\endaddress
\email andrew\@sophie.math.uga.edu \endemail
\address (Poonen) Department of Mathematics, The University of California,
Berkeley, California 94720-3840, USA \endaddress
\email poonen\@math.berkeley.edu \endemail
\address{(Soundararajan) Department of Mathematics, 
Princeton University, Princeton, New Jersey 08544, USA} \endaddress
\email{skannan{\@}math.princeton.edu} \endemail
\thanks The first three authors are all supported, in part, by grants from the
N.S.F. The first and fourth authors are supported by the American Institute of 
Mathematics. The second author is a Presidential Faculty Fellow.
The third author is an Alfred P. Sloan Research Fellow and a Packard Fellow.
 \endthanks
\leftheadtext{Conrey et. al.}

\endtopmatter

\document

\head 1. Introduction \endhead

\noindent Dirichlet noted that, from the formula
$$
\Gamma(s) = n^s\int^{\infty}_0 x^{s-1} e^{-nx} dx = n^s\int^1_0(-\log
t)^{s-1}t^{n-1}dt ,
$$ 
we may obtain the identity
$$
\align
\Gamma(s)L(s, (\tfrac{\cdot}{p}))
&= \Gamma(s) \sum_{n\ge 1}\frac{\L{n}{p}}{n^s} =
\int^1_0(-\log t)^{s-1}\sum_{n\ge 1}\L{n}{p} t^{n-1} dt\\ 
&= \int^1_0 \frac{(-\log t)^{s-1}}{t} \frac{f_p(t)}{1 - t^p} dt.
\tag{1.1}\\
\endalign
$$
Here $\L{\cdot}{p}$  is the Legendre symbol and 
$$
f_p(t) := \sum^{p-1}_{a=0} \L{a}{p} t^a. \tag{1.2}
$$
Equation (1.1) allowed Dirichlet to define $L(s, \L{\cdot}{p})$ 
as a regular function for all complex $s$.  Fekete observed that if 
$f_p(t)$ has no real zeros $t$ with $0<t<1$, then 
$L(s, \L{\cdot}{p})$ has no real zeros $s>0$;
and the $f_p(t)$ are thus now known as {\sl Fekete polynomials}.
Indeed, if $L(s, \L{\cdot}{p}) = 0$ then by (1.1) and the mean value theorem 
there is a $t$ in $(0,1)$ with $\frac{(-\log t)^{s-1}}{t} \frac{f_p(t)}{1-t^p}
=0$, and so $f_p(t)=0$ here.

Among small primes $p$, there are only a few for which the 
Fekete polynomial $f_p(t)$ has a real zero $t$ in the range $0<t<1$. 
In fact, we may verify computationally that there are 
just $23$ primes up to $1000$ for which $f_p$ has a zero in $(0,1)$.  
This implies that there are no  positive real zeros of 
$L(s, \L{\cdot}{p})$ for most such primes
$p$, and in particular no {\sl Siegel zeros} (that is, 
real zeros ``especially close to $1$'').
It is interesting to note that for those primes
$p\equiv 3 \mod 4$ for which $f_p(t)$ does have a zero in $(0,1)$, 
the class number of $Q(\sqrt{-p})$ is surprisingly small (for example 
$p=43, 67, 163, \dots$). Unfortunately this trend does not
persist: Indeed Baker and Montgomery [1] proved that $f_p(t)$
has a large number of zeros in $(0,1)$ for almost all primes $p$ (that
is, the number of such zeros $\to \infty$ as $p\to \infty$, and it seems
likely that there are, in fact, $\asymp \log\log p$ such zeros).

In this paper we shall study the complex zeros of $f_p(t)$.  Using
zero locating software one finds that, for primes $p$ up to $1000$,
about half of the zeros lie on the unit circle; leading one to expect this
to be the general phenomenon.  It turns out to be fairly easy to prove
that {\sl at least} half of the zeros of $f_p(t)$ are on the unit circle
(that is $|t|=1$):  First note that
$$
F_p(z) := z^{-p/2} f_p(z) = \sum^{(p-1)/2}_{a=1}
\L{a}{p}\left(z^{a-p/2} + \L{-1}{p} z^{p/2-a}\right)
$$
by combining the  $a$  and  $p-a$  terms\footnote{Here $z = e^{2i\pi t}$
 with $0\leq t<1$, so that there is no ambiguity in the meaning
of $z^{-p/2}$.}.  Taking  $z = e^{2i\pi t}$  we have
$$
F_p\left(e^{2i\pi t}\right) = 
\cases
2 \sum^{(p-1)/2}_{a=1} \left(\frac{a}{p}\right) \cos((2a-p)\pi t)
&\text{if  $p \equiv 1 \mod{4}$}\\
2i \sum^{(p-1)/2}_{a=1} \left(\frac{a}{p}\right) \sin((2a - p)\pi
t) &\text{if  $p \equiv 3 \mod{4}$}.
\endcases
\tag{1.3}
$$
Define $H_p(t)=F_p(e^{2i\pi t})$ if $p\equiv 1 \pmod 4$, and
$H_p(t)=-iF_p(e^{2i\pi t})$ if $p\equiv 3 \pmod 4$. By (1.3)
we see that $H_p(t)$ is a periodic, continuous real-valued function when 
$t$ is real.

Now if  $\zeta_p = e^{2i\pi/p}$  then, for all  $k$  not divisible
by  $p$,  $f_p(\zeta^k_p)$  is a {\sl Gauss sum} and has absolute
value  $\sqrt{p}$ (see section 2 of [2]);  therefore $|F_p(\zeta^k_p)| = \sqrt{p}$. 
Moreover
$$
F_p(\zeta^k_p) = (\zeta^k_p)^{-p/2} \sum^{p-1}_{a=1}
\L{a}{p} \zeta^{ak}_p 
= (-1)^k \L{k}{p} \sum^{p-1}_{a=1}
\L{ak}{p} \zeta^{ak}_p =(-1)^k\L{k}{p}F_p(\zeta_p) .
$$
Therefore if  $\L{k}{p} = \L{k+1}{p}$ 
then  $H_p(k/p)$  and  $H_p((k+1)/p)$  have different
signs.  Since  $H_p(t)$  is real-valued and continuous, it must 
have a zero in-between $k/p$ and $(k+1)/p$, 
by the intermediate value theorem.  
Thus the number of zeros of  $H_p(t)$ in $[0,1)$ (and so of $F_p(z)$  
on the unit circle) is 
$$
\ge \#\biggl\{k : 1 \le k \le p - 2\ \text{\rm  and}\ \L{k}{p}= 
\L{k+1}{p}\biggr \} =
\frac{p-3}{2}, 
$$
as we shall see in Lemma 2.

Other than possible zeros at  $z = -1$  and at  $z = 1$,  this
accounts for all the zeros on the unit circle for each prime  $p < 500$.  So
the question is, is this all, for all  $p$?  The answer is ``no'' and indeed
one finds more zeros when $p=661$. In general one has the following:

\proclaim{Theorem 1}  There exists a constant  $\kappa_0$, $1 > \kappa_0 >
\frac{1}{2}$  such that
$$
\#\{z : |z| = 1 \text{\rm \ and } f_p(z) = 0\} \sim \kappa_0p \text{\rm\
as } p \rightarrow \infty .
$$
\endproclaim

We determine  $\kappa_0$  in terms of another constant  $\kappa_1$  defined
as follows:

\proclaim{Theorem 2}  Let  $\Cal F_J$  be the set of rational functions 
$$
g(x) = \frac{1}{x} + \frac{1}{1-x} + \sum\Sb |j| < J \\ j\ne 0,-1
\endSb \frac{\delta_j}{x+j} 
$$
where we allow  each  $\delta_j$ to take value $+1$ or $-1$.  
There exists a constant $\kappa_1, \frac{1}{2} > \kappa_1 > 0$, such that
$$
\#\{g \in \Cal F_J : g(x) = 0 \text{\rm \  for some } x \in (0,1)\}
\sim \kappa_1 \#\{g \in \Cal F_J\}
$$
as $J\to \infty$.  \endproclaim

The constants $\kappa_0$ and $\kappa_1$ are related as follows:

\proclaim{Theorem 1$\frac{1}{2}$}  
In fact $\kappa_0 = \frac{1}{2} + \kappa_1$.
\endproclaim

It is still an open question to determine the value of $\kappa_0$.
It is known that a ``random'' 
trigonometric polynomial of degree $p$
has $p/\sqrt{3}$ zeros in $[0,1)$ (see [7]), 
so one might guess that 
$\kappa_0=1/\sqrt{3}\approx 0.5773\dots$. However this is not the case.
We will show
$$
0.500813 > \kappa_0 > 0.500668.  
$$
While it is theoretically easy to find the value of $\kappa_0$, we 
do not know a good practical way of achieving this.  

As well as determining precisely the proportion, $\kappa_0$, of the zeros of
$f_p(t)$ which lie on the unit circle, we would also like to understand the
distribution of the set of zeros in the complex plane. There are several 
easy remarks to make:\ By (1.2) we have 
$$
t^pf_p(1/t) = \left(\frac{-1}{p}\right)f_p(t)
$$
and so the zeros of  $f_p(t)$, other than  $t = 0$, are symmetric
about the unit circle (i.e. they come in pairs other than at 
$t=0, \pm 1$).  We also note that, for  $|t| > 1$,  
$$
|f_p(t)/t^{p-1}| = \left| \sum^{p-1}_{a=0} \left(\frac{a}{p}\right)
\frac{1}{t^{p-1-a}} \right|  \ge 1 - \sum^{p-2}_{a=0}
\frac{1}{|t|^{p-1-a}} > 1 - \frac{1}{|t|-1} .
$$
However if  $|t| \ge 2$  then  $1 - 1/(|t|-1) \ge 0$, and so 
$f_p(t)$  has no zeros in  $|t| \ge 2$.  By symmetry it has no
zeros in  $|t| \le \frac{1}{2}$ except 0.  Thus

\proclaim{Proposition 1}  The zeros of  $f_p(t)$, other than at 
 $0, 1$  and  $-1$  come in pairs  $\alpha,1/\alpha$.  Moreover,
other than  $0$, they all lie in the annulus  $\{r \in \Bbb C :
\frac{1}{2} < |r| < 2\}$.
\endproclaim

As for the distribution of the arguments of the roots of $f_p(t)$ we can use 
a beautiful result of Erd\H os and Tur\' an (Theorem 1 of [3]), which immediately implies that, for any $0\leq \alpha < \beta <1$,
$$
\# \{ \tau\in {\Bbb C}:\ f_p(\tau)=0, \ \alpha < \arg(\tau)/2\pi < \beta\}
= (\beta-\alpha) p + O(\sqrt{p\log p}). \tag{1.4}
$$
The arguments above, and those used in proving Theorems 1 and 2, focus on determining which arcs $(\zeta_p^K,\zeta_p^{K+1})$ of the unit circle
contain a zero of $f_p(t)$.  Evidently (1.4) cannot be used so precisely.
However we can show that there are zeros of $f_p(t)$ near to such an
arc, so long as $f_p(t)$ gets ``small'' on that arc.

\proclaim{Theorem 3}  Suppose that $\epsilon>0$ is a sufficiently small
constant. If $p$ is a sufficiently large prime and $K$ an integer such that 
there exists a value of $t$ on the unit circle in the arc from 
$\zeta_p^K$ to $\zeta_p^{K+1}$ with $|f_p(t)|<\epsilon \sqrt{p}$, then there
exists $\tau=r\zeta_p^{K+\theta}$ with $f_p(\tau)=0$ where 
$0<\theta <1$ and $1-\epsilon^{1/3}/p < r \leq 1$.
\endproclaim

\remark{Remark} Applying Proposition 1 we also have $f_p((1/r)\zeta_p^{K+\theta})=0$. \endremark

As we have already discussed, Gauss sums 
 $\sum^{p-1}_{a=1} \L{a}{p} \zeta^{ak}_p$ (and many
generalizations) have the surprising property that they have absolute value
exactly equal to $\sqrt{p}$. It is, we think, of interest to ask
what happens when we replace the primitive $p$th root of unity $\zeta^{k}_p$
in the expression for a Gauss sum above, by some primitive $2p$th root of 
unity. These may be written as $\zeta^{k+1/2}_p$ or $\zeta^{2k+1}_{2p}$, or
$-\zeta^{k}_p$; so we must consider the values of $f_p(-\zeta^{k}_p)$.
Do these all take on the same absolute value ? The answer we now see
is ``no'', as we evaluate the distribution of these absolute values:

\proclaim{Theorem 4}  For any fixed real number  $\rho$
$$
\#\biggl\{k : 1 \le k \le p \text{\rm\ such that }
H_p\left(\frac{k+1/2}{p}\right) < \rho\sqrt{p}\biggr\} \sim c_{\rho} p
$$
as $p\to \infty$ where
$$
c_{\rho} = \frac{1}{2} + \frac{1}{\pi} \int^{\infty}_{x=0}
\sin(\rho \pi x) \prod\Sb n \geq 1 \\ n \text{\rm \ odd} \endSb
\cos^2\left(\frac{2x}{n}\right) \frac{dx}{x} .
$$
Moreover $c_{-\rho}$ and $1-c_\rho=\exp(- \exp( \pi \rho/2 + O(1)))$ for 
positive $\rho$.
\endproclaim

After proving this in section 6, we indicate how our proof may be
modified to establish several related results. First, to show that 
$\max_{|z|=1} |f_p(z)| \gg \sqrt{p} \log \log p$, so re-establishing
a result of Montgomery [5]. Second to understand the distribution
of the values of the Fekete polynomial at $(p-1)$st roots of unity.

\noindent{\smc Acknowledgements}:  {\sl We thank Jeff Lagarias for
facilitating this joint endeavour, Peter Borwein, Neil Dummigan, Hugh Montgomery,  Pieter Moree, Mike Mossinghoff, Bob Vaughan and Trevor Wooley
for some helpful remarks, and the referee for a very careful reading of the
paper.}

\head 2.  First results \endhead

\noindent Let $\chi$ be any character $\pmod p$ 
and let $k$ be an integer not divisible by $p$.  Note that 
$$
\sum^{p-1}_{a=1} \chi(a) \zeta^{ak}_p = \bar\chi (k) \sum_{a=1}^{p-1} 
\chi(ak) \zeta_p^{ak} 
= \bar\chi(k) \sum^{p-1}_{b=1} \chi(b) \zeta^{b}_p. \tag{2.1}
$$
In particular we see that $ f_p(\zeta^k_p) = \fracwithdelims(){k}{p} 
f_p(\zeta_p)$, whereas in contrast $f_p(1) = 0$.  Recall that for a 
non-principal  character $\chi\pmod p$, the Gauss sum $\tau(\chi)$ is
$\sum_{a=1}^{p-1} \chi(a) \zeta_p^a$.  Thus $f_p(\zeta_p)$ is the 
Gauss sum $\tau(\L{\cdot}{p})$.  It is easy to determine 
the magnitude of $|f_p(\zeta_p)|$:  Note that 
$$
\align
(p-1)f_p(\zeta_p)^2 &= \sum^{p-1}_{k=0} f_p(\zeta^k_p)^2 =
\sum^{p-1}_{k=0} \sum^{p-1}_{a,b=0} \L{ab}{p}\zeta^{(a+b)k}_{p} \\
&= \sum^{p-1}_{a,b=1} \L{ab}{p}\sum^{p-1}_{k=0} 
\zeta^{(a+b)k}_p = p \sum\Sb a=1 \\ b =
p-a\endSb^{p-1} \L{ab}{p}  =p\L{-1}{p}(p-1).
\endalign
$$
Hence we have $f_p(\zeta_p)^2 = \L{-1}{p} p$, and so $|f_p(\zeta_p)|=\sqrt{p}$.
Gauss showed more and determined that
$$
f_p(\zeta_p) = 
\cases 
\sqrt{p} &\text{if  } p\equiv 1\pmod 4,\\
i\sqrt{p} &\text{if  } p\equiv 3\pmod 4.\\
\endcases
$$

Since $f_p(\zeta^k_p) =\L{k}{p}f_p(\zeta_p)$, for  $1 \le k \le p-1$, and 
$f_p(1) = 0$, we get by Lagrangian interpolation 
$$
f_p(z) = \sum^{p-1}_{k=0} f_p(\zeta^k_p) \prod\Sb j=0 \\ j
\ne k\endSb^{p-1} \biggl(\frac{z-\zeta^j_p}{\zeta^k_p - \zeta^j_p}\biggr).
$$
Note that 
$$
\prod\Sb j=0\\ j\neq k\endSb^{p-1} 
(z-\zeta_p^j) = \frac{z^p-1}{z-\zeta_p^k}, \qquad \text{and that} 
\qquad
\prod\Sb j=0\\ j\neq k\endSb^{p-1} (\zeta_p^k-\zeta_p^j) = \zeta_p^{k(p-1)}
\prod_{j=1}^{p-1} (1-\zeta_p^j) = p \zeta_p^{-k}.
$$ 
Hence 
$$
\frac{p}{f_p(\zeta_p)} \frac{ f_p(z)}{z^p - 1} 
= \frac{p}{f_p(\zeta_p)} \frac{z^{-\frac p2} f_p(z)}{z^{\frac p2}-z^{-\frac 
p2}} = \sum^{p-1}_{k=1} \L{k}{p}\frac{\zeta^k_p}{z - \zeta^k_p}.
\tag{2.2}
$$
If $|z|=1$ then note that $z^{\frac p2}-z^{-\frac p2} \in i{\Bbb R}$,
and from (1.3) and $f_p(\zeta_p)^2=\L{-1}{p}p$ we have $z^{-\frac p2} f_p(z)
/f_p(\zeta_p) \in {\Bbb R}$.  Thus the right side of (2.2) $\in i{\Bbb R}$
for all $|z|=1$.  To facilitate studying $f_p(z)$ as $z$ goes 
around the unit circle from $\zeta^{K}_p$ to $\zeta^{K+1}_p$, we write 
$z = \zeta^{K + x}_p = \zeta^{K}_p e^{2 i \pi x/p}$ and then let 
$$
g_{p, K}(x) := i\L{K}{p}\frac{p}{f_p(\zeta_p)} \frac{f_p(z)}{z^p - 1} 
\biggr|_{z = \zeta^{K +x}_p} 
= i\L{K}{p} \sum^{K+(\frac{p-1}{2})}_{k = K -(\frac{p-1}{2})}
\L{k}p \frac{1}{\zeta^{K - k + x}_p - 1}. \tag{2.3}
$$
Thus $g_{p,K}(x)$ is a real valued function of $x\in [0,1]$.  

\proclaim{Proposition 2} If $0\le K\le p-1$ is an integer with $\L{K}{p}
=\L{K+1}p$ then $g_{p,K}(x)$ has exactly one zero in $(0,1)$.  
Equivalently, $f_p(z)$ has exactly one zero on the arc of the 
unit circle from $\zeta_p^K$ to $\zeta_p^{K+1}$.  If $\L{K}p =
-\L{K+1}p$ then $g_{p,K}$ has either no zeros, or exactly two 
zeros in $(0,1)$.  Equivalently, $f_p(z)$ has exactly $0$ or $2$ 
zeros on the arc from $\zeta_p^K$ to $\zeta_p^{K+1}$.
\endproclaim

\remark{Remark}  In the above Proposition, and henceforth, we count zeros 
with multiplicity.
\endremark

Before proving the Proposition, we evaluate $\sum_{k=1}^{p-1} 
\frac{1}{\sin^2 (\frac{\pi k}{p})}$.

\proclaim{Lemma 1} For all integers $p\ge 2$,
$$
\sum_{k=1}^{p-1} \frac{1}{\sin^2 (\frac{\pi k}{p})} 
=\frac{p^2-1}{3}.
$$
\endproclaim
\demo{Proof}  Put $A(z) = \prod_{k=1}^{p-1} (z-\zeta_p^k)$.  
Logarithmic differentiation shows that 
$$
\left\{ z \L{A^{\prime}(z)}{A(z)}^{\prime}+ \frac{A'(z)}{A(z)}\right\} \ \biggl|_{z=1} 
= -\sum_{k=1}^{p-1} \frac{\zeta_p^k}{(1-\zeta_p^k)^2} = \frac{1}{4} 
\sum_{k=1}^{p-1} \frac{1}{\sin^2 (\frac {\pi k}{p})}.
$$
However, $A(z) = \frac{z^p-1}{z-1} = z^{p-1}+z^{p-2}+\ldots +1$ 
and using this to evaluate the left side above, we get the lemma.
\enddemo

\demo{Proof of Proposition 2} Note that with $g=g_{p, K}$, we have 
$\lim_{x\to 0^+} g(x) =  \infty$, and 
$\lim_{x\to 1^-} g(x) = - \L{K}{p}\L{K+1}{p}\infty$.
Further observe that 
$$
\align
g'(x) &= \frac{2\pi}{p}  \L{K}{p} 
\ \sum_{|k - K| < \frac{p}{2}} \ \L{k}{p} \frac{\zeta_p^{K-k+x}}
{(\zeta^{K-k+x}_p -1)^2}\\
&= - \frac{\pi}{2p} \L{K}{p} \ \sum_{|k-K| < \frac{p}{2}}\  \L{k}{p} 
\frac{1}{\sin^2 (\frac{\pi}{p}(K-k+x))}.\\
\endalign
$$

If  $\L{K}{p} = \L{K+1}{p}$ then, by Lemma 1,
$$
\align
|g'(x)| &\ge \frac{\pi}{2p} \biggl(
\frac{1}{\sin^2 (\frac{\pi}{p}x)} + \frac{1}{\sin^2 (\frac{\pi}{p}(1-x))} 
- \sum\Sb j\ne 0,1 \\
|j|<p/2\endSb \ \frac{1}{\sin^2(\frac{\pi}{p}(x-j))}\biggr)\\
& \geq \frac{\pi}{2p} \biggl( 
\frac{2}{\sin^2 (\frac{\pi}{2p})} - \frac{p^2-1}{3} \biggr)
> 0, \tag{2.4}\\
\endalign
$$
since the sum of the first two terms is minimized when $x=\frac 12$.
Hence $g'(x) \ne 0$  for all  $x \in (0,1)$, so that $g$  is
monotone decreasing in $[0,1]$ going from  $\infty$ to $-\infty$.  Thus
$g$ has exactly one zero in this interval.

Moreover
$$
g^{\prime\prime}(x) =
\frac{\pi^2}{p^2} \L{K}{p}
\ \sum_{|k-K|<p/2} \L{k}{p} \frac{\cos (\frac{\pi}{p}(K-k+x))}
{\sin^3  (\frac{\pi}{p}(K-k+x))}.  
$$
Now if  $\L{K}{p} = -\L {K+1}{p}$ then
$$
g^{\prime\prime}(x) \ge \frac{\pi^2}{p^2} \biggl(
\frac{\cos (\frac{\pi}{p}x)}{\sin^3 (\frac{\pi}{p} x)} 
+\frac{\cos (\frac{\pi}{p}(1-x))}{\sin^3 (\frac{\pi}{p}(1-x))} 
- \sum\Sb |j| < p/2 \\ j \ne 0,1\endSb 
\frac{\cos (\frac{\pi}{p} (j-x))}{|\sin (\frac{\pi}{p}(j-x))|^3}\biggr).
$$
Let $\mu$ be the minimum of $\cot (\frac {\pi}{p} t)$ over $t=x$, $1-x$.
Since $\cot t$ decreases rapidly as $t$ goes from $0$ to $\frac {\pi}{2}$ 
we see that the above is 
$$ 
\geq \frac{\pi^2}{p^2}\mu  \biggl(
\frac{1}{\sin^2 (\frac{\pi}{p}x )} 
+ \frac{1}{\sin^2 (\frac{\pi}{p}(1-x))} 
- \sum\Sb j\ne 0,1 \\ |j|<p/2\endSb \frac{1}{\sin^2 (\frac{\pi}{p}(x-j))}
\biggr) > 0,
$$
as in (2.4).  Thus  $g'(x)$  is monotone increasing in  $(0,1)$ going 
from $- \infty$  to  $+ \infty$.  Thus there is a unique $x_0$ in $(0,1)$ 
with $g'(x_0)=0$, and the minimum value of $g(x)$ is attained at $x_0$.  
Plainly $g$ has $0$ or $2$ zeros depending on whether $g(x_0)>0$, or 
$g(x_0)\le 0$.  This proves the proposition.
\enddemo

>From Proposition 2 we know that $f_p(z)$ has at least as many zeros on $|z|=1$,
as there are values $1\le K\le p-1$ with $\L{K}{p}=\L{K+1}{p}$.  
We next determine the number of such values $K$.

\proclaim{Lemma 2 (Gauss)} For any non-principal character $\chi \pmod p$, 
we have
$$
\sum^{p-1}_{b=1} \chi(b) \bar\chi(b+k) = 
\cases 
p-1 &\text{if  } p | k \\ 
-1  &\text{if  } p \nmid k.\\
\endcases
\tag{2.5}
$$
Hence 
$$
\#\left \{ b \pmod{p} :\L{b}{p} = \L{b+1}{p}\right \} =\frac{p-3}{2},
$$
and 
$$
\#\left \{ b \pmod{p} : \L{b}{p}= -\L{b+1}{p}\right \} = \frac{p-1}{2}.
$$
\endproclaim

\demo{Proof}  If  $p|k$  then the right side of (2.5) is 
$\sum_{b=1}^{p-1}|\chi(b)|^2 =p-1$.  Suppose now that $p\nmid k$, and let 
$c=(b+k)/b = 1+k/b$.  As $b$ runs over the non-zero residue classes $\pmod q$,
note that $c$ runs over all residue classes except the residue class 
$1 \pmod p$.  Hence the right side of (2.5) is 
$$
\sum\Sb c\pmod p \\ c\not\equiv 1 \pmod p\endSb \bar\chi(c) = -1,
$$
as desired.
\enddemo

If $\L{K}{p}=-\L{K+1}{p}$ then we need to determine (in the notation 
of the proof of Proposition 2) whether $g(x_0)>0$ or $\le 0$.  
This depends heavily on the values of $\L{k}{p}$ for $k$ neighbouring 
$K$.  The following Lemma shows that these neighbouring values behave like 
independent random variables.

\proclaim{Lemma 3 (Weil)}  Fix integer $J$, and then the numbers
 $\delta_j \in \{-1,1\}$  for each  $j$ with $|j| < J$.  We have, uniformly, 
$$
\#\biggl\{x \pmod{p} : \L{x-j}{p} = \delta_j\ 
\text{for all } \ |j| < J\biggr\} = \frac{p}{2^{2J-1}} + O(J\sqrt{p}).
$$
\endproclaim
\demo{Proof}  The above equals
$$
\align
\sum^p_{x=1} \frac{1}{2^{2J-1}} &\prod_{|j| < J }
\left(1 + \delta_j\left(\frac{x-j}{p}\right)\right) +
O(J) \\
&= \frac{p}{2^{2J-1}} + O\biggl(\frac{1}{2^{2J-1}} \sum\Sb
S\subseteq \{ |j|<J\}\\ S \ne \emptyset\endSb \ \ \sum^p_{x=1}
\L{\prod_{j\in S} (x-j)}{p} + J \biggr).
\endalign
$$
By Weil's Theorem [8], if  $f(x)$  is a squarefree polynomial 
$\pmod{p}$  then
$$
\left|\sum^p_{x=1} \left(\frac{f(x)}{p}\right)\right| \ll
(\text{degree } f)\sqrt{p}.
$$
Hence the above is
$$
= \frac{p}{2^{2J-1}} + O\left(\frac{\sqrt{p}}{2^{2J-1}}
\sum^{2J-1}_{m=1} \pmatrix 2J-1 \\ m\endpmatrix m + J\right),
$$
and the result follows.
\enddemo

We conclude this section by determining the order of the zeros of 
$f_p(z)$ at $\pm 1$.  In fact we shall determine the number of zeros of 
$f_p(z)$ on the arcs $\zeta_p^{\frac{p-1}{2}}$ to $\zeta_p^{\frac{p+1}{2}}$
(which contains $-1$), and $\zeta_p^{-1}$ to $\zeta_p$ 
(which contains $1$).

\proclaim{Lemma 4}  If $p\equiv 1\pmod 4$ then $f_p(z)$ has only a
simple zero at $z=-1$, on the arc from
$\zeta_{p}^{\frac{p-1}{2}}$ to $\zeta_{p}^{\frac {p+1}{2}}$, and
$f_p(z)$ has only a double zero at $z=1$, on the arc from 
$\zeta_p^{-1}$ to $\zeta_p$.  
If $p\equiv 3 \pmod 4$ then there are no zeros of $f_p(z)$ on the
arc from $\zeta_{p}^{\frac{p-1}{2}}$ to $\zeta_{p}^{\frac {p+1}{2}}$,
and $f_p(z)$ has only a simple zero at $z=1$ on the arc from $\zeta_p^{-1}$ 
to $\zeta_p$.   
\endproclaim
\demo{Proof}  We make free use of the fact that $\L{-1}{p}= 1$, or $-1$ 
depending on whether $p\equiv 1\pmod 4$, or $3\pmod 4$.  
Let's begin with the arc from $\zeta_p^{\frac{p-1}{2}}$ to $\zeta_p^{\frac{p+1}
2}$.  We take $K=\frac{p-1}{2}$ in Proposition 2.  Note that 
$\L{K}{p}=\L{K+1}{p}$ if $p\equiv 1\pmod 4$, and $\L{K}{p}=-\L{K+1}{p}$
if $p\equiv 3 \pmod 4$.  In the first case, Proposition 2 tells
us that there's exactly one (simple) zero on this arc.  Since 
$f_p(-1)=\sum_{a=1}^{p-1} (-1)^a \L{a}{p} = \frac 12 \sum_{a=1}^{p-1} 
(-1)^a (\L{a}{p}-\L{p-a}{p}) = 0$ for $p\equiv 1\pmod 4$, this 
simple zero is at $-1$.  Now suppose $p\equiv 3 \pmod 4$.  By 
Proposition 2, we know that there are $0$ or $2$ zeros on this arc, 
depending on whether $\min_x g_{p,K}(x) >0$ or not.  We now show 
that this minimum is attained at $x=\frac 12$, and the minimum 
value is positive.  Putting $j=K-k$ in (2.3) we have
$$
\align
g_{p,K}(x)&= i\L{K}{p} \sum_{|j|\le \frac{p-1}{2}} 
\L{K-j}{p} \frac{1}{\zeta_p^{j+x}-1} \\
&= i\L{K}{p} \sum_{j=0}^{\frac{p-1}{2}} \L{K-j}{p} 
\biggl(\frac{1}{\zeta_p^{j+x}-1} -\frac{1}{\zeta_p^{-j-1+x}-1}\biggr),\\
\endalign
$$
since $K+j+1\equiv -(K-j) \pmod p$. Evidently 
$g_{p,K}(1-x)=\overline{g_{p,K}(x)}$, so $g_{p,K}(1-x)=g_{p,K}(x)$
since $g_{p,K}(x)$ is real-valued. However we see that the minimum
of $g_{p,K}(x)$ is obtained at a unique point in $(0,1)$, so that 
must be at  $x=\frac 12$.  Now 
$$
f_p(-1) = \sum^{p-1}_{a=1} (-1)^a \L{a}{p} =
\sum\Sb a=1\\ a \ \text{even}\endSb^{p-1} \L{a}{p}- 
 \sum\Sb b=1\\b \ \text{even}\endSb^{p-1} \L{p-b}{p}
$$
where  $a = p - b$  is odd in the second sum,
$$
= 2 \sum^{(p-1)/2}_{d=1} \left(\frac{2d}{p}\right) =
2\left(\frac{2}{p}\right) \sum^{(p-1)/2}_{d=1}
\left(\frac{d}{p}\right) = 2\left( 2\left(\frac{2}{p}\right) -1\right) 
h(-p),
$$
where $h(-p)$ is the class number of ${\Bbb Q}(\sqrt{-p})$ 
(see section 2 of [2]).  By (2.3), and since $f_p(\zeta_p)=i\sqrt{p}$ 
by Gauss, we have 
$$
\align
g_{p,K}(\tfrac 12)&= - \L{K}{p}\frac{\sqrt{p}}{2} f_p(-1) = \sqrt{p} 
\biggl(-2 \L{2K}{p} + \L{K}{p}\biggr) h(-p)\\ &= 
\sqrt{p} \biggl(2 +\L{K}{p}\biggr) h(-p) >0. \\
\endalign
$$
This shows that $f_p(z)$ has no zeros on the arc from 
$\zeta_p^{\frac{p-1}{2}}$ to $\zeta_p^{\frac{p+1}{2}}$ when 
$p\equiv 3\pmod 4$.

Now let's consider the arc from $\zeta_p^{-1}$ to $\zeta_p$.  
Take $K=p-1$, and consider $g_{p,K}(x)$ as defined in (2.3).  
Usually $g_{p,K}(x)$ would have a discontinuity at $1$, 
but here since $\L{K+1}{p}=\L{0}{p}=0$ we do not have this problem.
Thus $g_{p,K}$ is a continuous function on $(0,2)$, and we 
may study $f_p(z)$ on the arc from $\zeta_p^{-1}$ to $\zeta_p$ 
by studying $g_{p,K}(x)$ on $(0,2)$. 
Note that for any  $p$,  $f_p(1) = \sum^{p-1}_{a=1}
\L{a}{p} = 0$, so that there is at least a simple 
zero at $z=1$. Also $f_p^\prime(1)=-i (-1/p) f_p(\zeta_p) g_{p,p-1}(1)$
by (2.3). Since $f_p(z)=(-1/p)z^pf_p(\overline{z})$, we deduce that 
$g_{p,p-1}(x)=-(-1/p)g_{p,p-1}(2-x)$.

If $p \equiv 1\pmod 4$ then $g_{p,p-1}(1)=0$ and so $f_p^\prime(1)=0$.
Now, as in the proof of (2.4), the first part of the proof of Proposition 2,
we have $|g^{\prime}_{p,K}(x)| >0$ for all $x\in(0,2)$.  Therefore
$g$ has only a simple zero at $x=1$, and thus $f_p$ has a double zero at $1$.

If $p\equiv 3\pmod 4$ then, as in the second part of the proof of Proposition 2,
$|g^{\prime\prime}_{p,K}(x)| >0$ for $x\in(0,2)$.  Thus there is a unique
minimum of $g_{p,K}(x)$ on $(0,2)$, but since $g_{p,p-1}(x)=g_{p,p-1}(2-x)$
this must be attained at $x=1$. However, by (2.3), and as $f_p(\zeta_p)
=i\sqrt{p}$ by Gauss,
$$
g_{p,K}(1) = -\frac{f_p^{'}(1)}{\sqrt{p}}= - \frac{1}{\sqrt{p}}
\sum_{a=1}^{p-1} a \L{a}{p} = 
\sqrt{p}h(-p) > 0,
$$
(see [2], section 2), and so $g_{p,K}(x)>0$ and thus has no zeros in $(0,2)$.
Therefore $f_p$ has only a simple zero at $z=1$ on this arc.  
\enddemo

\head 3.  Functions with random coefficients \endhead

\noindent If $g \in \Cal F_J$ then, for any  $x \in (0,1)$, we
have
$$
\align
\frac{1}{2} g^{\prime\prime}(x) &= \frac{1}{x^3} +
\frac{1}{(1-x)^3} + \sum\Sb |j|< J \\ j\ne 0,-1\endSb
\frac{\delta_j}{(x+j)^3} \\
&\ge \frac{1}{x^3} + \frac{1}{(1-x)^3} - \sum\Sb |j|< J \\ j\ne
0,-1\endSb \frac{1}{(x+j)^3} > 2 \frac{1}{(\frac 12)^3} - 2\zeta(3) > 0.
\tag{3.1}\\
\endalign
$$
Since  $\lim\limits_{t\rightarrow 0^+} g'(t) = -\infty$  and 
$\lim\limits_{t\rightarrow 1^-} g'(t) = \infty$  we deduce that 
$g'(x)$  has exactly one zero in  $(0,1)$, call it  $x_0$.  Note that
$g(x)$ attains its minimum value at $x_0$.
If  $0 \le t < \frac {1}{\pi}$  then
$$
-g'(t) \ge \frac{1}{t^2} - 2\left(\frac{1}{(1/2)^2} +
\frac{1}{(3/2)^2} + \frac{1}{(5/2)^2} + \dots \right) =
\frac{1}{t^2} - \pi^2 > 0. 
$$
Similarly if  $1 - \frac {1}{\pi} < t \le 1$ then $g'(t)>0$.  Thus
$$
x_0 \in [\tfrac {1}{\pi},1 - \tfrac{1}{\pi}]. \tag{3.2}
$$
We now show that few $g$ are small in absolute value, at their minimum $x_0$.

\proclaim{Proposition 3} We have $|g(x_0)|> J^{-\frac 14}$ for almost
all $g \in \Cal F_J$, where $g'(x_0)=0$, uniformly as $J\to \infty$.
\endproclaim
\demo{Proof}  Consider the subset $S$ of  $\Cal F_J$  
with all the $\delta_j$  fixed given values,
except when  $j \in [I,I+I^{\frac 12}]$  where  $I =J^{\frac 14}$.
Let  $f \in S$  with  $\delta_j = -1$  for all  $j \in [I,I
+ I^{\frac 12}]$.  Suppose that  $f'(x_1)=0$ and let 
$$
\gamma = \sum\Sb |j| < J \\ j \notin [I,I+I^{\frac 12}]\endSb
\frac{\delta_j}{x_1+j}
$$
where  $\delta_0 = 1$, $\delta_{-1} = -1$.  Let $g$ be any element of
$S$ with $g'(x_0)=0$.

By (3.1) note that 
$$
\align 
|x_1 - x_0| &\ll \biggl|\int_{x_0}^{x_1} f^{\prime\prime}(t)dt \biggr| 
= |f^{\prime}(x_0)-f^{\prime}(x_1)| = |f^{\prime}(x_0)| \\
&=|f^{\prime}(x_0)-g^{\prime}(x_0)| \le 2\sum_{j\in [I,I+I^{\frac12}]}
\frac{1}{(x_0+j)^2} \ll \frac 1I.\tag{3.3}
\\
\endalign
$$
Hence, keeping in mind $x_0$, $x_1\in [\frac 1\pi, 1-\frac 1\pi]$,
$$
\align
g(x_0) - \gamma &= \sum_{j\in [I,I+I^{\frac 12}]} 
\frac{\delta_j}{x_0 +
j} + O\biggl(\sum\Sb |j|< J\\ j \notin [I,I+I^{\frac 12}]\endSb
 \left|\frac{1}{x_0+j} - \frac{1}{x_1+j}\right|\biggr) \\
&= \frac{1}{I} \sum_{j\in[I,I+I^{\frac 12}]} \delta_j + O\biggl(
\sum_{j\in[I,I+I^{\frac 12}]} \left| \frac{1}{I} - \frac{1}{x_0+j}\right|
+|x_1-x_0| \biggr) \\
&= \frac{1}{I} \sum_{j\in[I,I+I^{\frac 12}]} \delta_j + 
O\biggl(\frac{1}{I}\biggr),
\\
\endalign
$$
since each $|1/I - 1/(x_0+j)| \ll 1/I^{\frac 32}$ and there are 
$I^{\frac 12}$ such 
terms.
Therefore if  $|g(x_0)| \leq \frac{1}{I}$ then
$$
\sum_{j\in[I,I+I^{\frac 12}]} \delta_j = -\gamma I + O(1).
\tag{3.4}
$$
Now, the $\delta_j$ are independent binomial random variables, so
the distribution of their sum tends towards the normal distribution. 
Therefore 
the maximum probability for (3.4) to occur happens when $\gamma = 0$; and so
(3.4) holds with probability $O(I^{-\frac 14})$, for any  $\gamma$,
 implying Proposition 3.
\enddemo

\head 4.  Proof of Theorem 2 \endhead

\noindent Suppose that  $g \in \Cal F_J$ and 
$f \in \Cal F_K$, with  $J < K$, 
such that  the $\delta_j$ are the same in each for  $|j| < J$.
Select $x_0, x_1\in (0,1)$ so that 
$g'(x_0)=0$ and  $f'(x_1) = 0$.  Now
$$
|f(x_1) - f(x_0)| \leq \sum_{|j|<K} \left|\frac{1}{x_1+j} -
\frac{1}{x_0+j}\right| \ll \sum_{|j|<K} \frac{|x_1 - x_0|}{j^2+1} \ll
|x_1 - x_0|,
$$
since  $x_0$, $x_1 \in [1/\pi, 1- 1/\pi]$.  Arguing exactly as in (3.3),
we see that $|x_0-x_1| \ll \frac{1}{J}$, and so we have 
$$
|f(x_1) - f(x_0)| \ll \frac{1}{J}. \tag{4.1}
$$

We next consider the mean-square of 
$$
|f(x_0) - g(x_0)| = \biggl|\sum_{J\le |j|<K} \frac{\delta_j}{x_0+j}\biggr|.
$$
To do so we will need to sum over all 
$\delta=\{ \delta_j\}_{J\le |j| < K}\in \Delta_{J,K}$, that is
the set of all possibilities with each $\delta_j=-1$ or $1$
(note that there are $2$ possible values for each $\delta_j$ so
the set $\Delta_{J,K}$ has $2^{2K-2J}$ elements). With this
notation, the mean square is
$$
\align
\frac{1}{2^{2K-2J}} \sum\Sb\delta\in \Delta_{J,K} \endSb
&\biggl|\sum_{J\le |j|<K} \frac{\delta_j}{x_0+j}\biggr|^2 
\\
&= \sum_{J\le |j_1|,|j_2| < K} \frac{1}{(x_0+j_1)(x_0 + j_2)}
\frac{1}{2^{2K-2J}}  \sum\Sb \delta\in \Delta_{J,K}\endSb \delta_{j_1}\delta_{j_2} \\
&= \sum_{J\le |j|<K} \frac{1}{(x_0+j)^2} \asymp \frac{1}{J}.\\ 
\endalign
$$
Thus if  $\psi_J \rightarrow \infty$  as  $J \rightarrow \infty$ 
then
$$
\biggl|\sum_{J\le |j| < K} \frac{\delta_j}{x_0+j}\biggr| <
\frac{\psi_J}{J^{\frac 12}}, \tag{4.2}
$$
for almost all choices of the  $\delta_j$.  

Combining (4.1) and (4.2), we see that for almost all choices of $\delta_j$
($J\le |j|<K$) we have
$$
|f(x_1) - g(x_0)| \le |f(x_1) - f(x_0)| + |f(x_0) - g(x_0)| <
\frac{2\psi_J}{J^{\frac 12}}. \tag{4.3}
$$
Taking $\Psi_J=J^{\frac 14}/2$, and combining this with Proposition 3 we
see that for almost all  $g \in \Cal F_J$, and almost all 
extensions $f$ of $g$ to ${\Cal F}_K$, $f(x_1)$  has the same
sign as  $g(x_0)$.  Summing up over all $g \in \Cal F_J$ we deduce that 
$\omega_K=\omega_J+o(1)$, where 
$$
\omega_J:=\frac{\#\{g \in \Cal F_J : g(x) = 0 \text{\rm \  for some } x \in (0,1)\}} {\#\{g \in \Cal F_J\}},
$$
and the ``$o(1)$'' term depends only on $J$. Therefore
$\lim_{J\to \infty} \omega_J$ exists, and equals $\kappa_1$ say.

Strong bounds on $\kappa_1$, which imply those in the statement of
Theorem 2, are given in Proposition 6 in section 8.

Theorem 2 follows.

\head 5.  Proofs of Theorems 1 and 1$\frac{1}{2}$ \endhead

\noindent Let $1\le K\le p-1$ be an integer.  If $\L{K}{p}=\L{K+1}{p}$ 
then by Proposition 2 there is exactly one zero of $f_p(z)$ on the 
arc from $\zeta_p^K$ to $\zeta_p^{K+1}$; by Lemma 2 this happens for 
$\sim \frac{p}{2}$ values of $K$.  Suppose now that 
$\L{K}{p}=-\L{K+1}{p}$ so that $f_p(z)$ has either $0$ or $2$ zeros on 
the arc from $\zeta_p^K$ to $\zeta_p^{K+1}$ depending on whether 
$\min_{x\in (0,1)} g_{p,K}(x)$ is positive or not. To decide this 
question we need the following proposition:

\proclaim{Proposition 4}  Suppose $J\le \sqrt{p}$, and $J\to \infty$ as 
$p\to \infty$.  For almost all $1\le K\le p-1$ we have 
$$
g_{p,K}(x) = \frac{p}{2\pi} \L{K}{p} \sum_{|j| < J} \L{K-j}{p} \frac{1}{j+x}
+ O\biggl(\frac p{J^{\frac 13}}\biggr),
$$ 
uniformly for all $x\in (0,1)$.
\endproclaim
\demo{Proof}  Note that for  $J \le |j| < \frac p2$,
$$
\biggl|\frac{1}{\zeta^{j+x}_p - 1} - \frac{1}{\zeta^j_p - 1}\biggr|
= \biggl| \frac{\zeta^x_p - 1}{(\zeta^{j+x}_p - 1)(\zeta^j_p -
1)}\biggr| \asymp \frac{px}{j(j+x)} \ll \frac{p}{j^2},
$$
and, for  $|j| < J$,
$$
\frac{1}{\zeta^{j+x}_p - 1} = \frac{p}{2i\pi} \frac{1}{(j+x)} + O(1).  
$$
Hence, putting $j = K-k$ in (2.3), we have
$$
\align
g_{p,K}(x) &= i \L{K}{p} \sum_{|j| < \frac p2}
\L{K-j}{p} \frac{1}{\zeta^{j+x}_p - 1}
\\
&= \frac{p}{2\pi } \L{K}{p} \sum_{|j|<J} \L{K-j}{p}
\frac{1}{j+x} +  i \L{K}{p}\sum_{J\le |j|<\frac p2} \L{K-j}{p} 
\frac{1}{\zeta^j_p - 1} 
+ O\left(J +\frac{p}{J}\right).\\
\endalign
$$

We now show that the mean-square of the second term above is small,
which proves the Proposition.  By Lemma 2,
$$
\align
\sum^p_{K=1} \biggl|\sum_{J\le |j|<\frac p2} &\L{K-j}{p}
\frac{1}{\zeta^j_p - 1}\biggr|^2 \\
&=\sum_{J\le |j_1|,\ |j_2|<\frac p2} \ \frac{1}{(\zeta^{j_1}_p -
1)(\zeta^{-j_2}_p - 1)} 
\sum^p_{K=1} \L{K-j_1}{p}\L{K-j_2}{p}\\
&= p \sum_{J\le |j|<\frac p2} \frac{1}{|\zeta^j_p - 1|^2} -
\biggl|\sum_{J\le |j| <\frac p2} \frac{1}{\zeta^j_p -
1}\biggr|^2 \\
&\ll p \sum_{J\le |j|<p/2} \left(\frac{p}{j}\right)^2 +
\biggl(\sum_{J\le |j|<p/2} \frac{p}{j}\biggr)^2 
\ll \frac{p^3}{J} + p^2\log^2 p.\\
\endalign
$$
This proves the Proposition.
\enddemo

By Proposition 4 we know that for almost all $K$ with $\L{K}{p}=-\L{K+1}{p}$
the minimum value of $\frac{2\pi}{p} g_{p,K}(x)$ 
equals the minimum of  $\L{K}{p} \sum_{|j|<J} \L{K-j}{p}\frac{1}{j+x} +
O(J^{-\frac 13})$.  For such $K$ the minimum value of $g_{p,K}(x)$ is 
non-positive if and only if  the minimum of  $\L{K}{p}\sum_{|j|<J} 
\L{K-j}{p}\frac{1}{j+x}$ is non-positive, unless
$$
\L{K}{p}\sum_{|j|<J} \L{K-j}{p}
\frac{1}{j+x}\ll \frac{1}{J^{\frac 13}}.  \tag{5.1}
$$

Now choose $J= [\frac{\log p}{10}]$.  Given any choice of 
$\delta_j\in \{ -1,1\}$,  $0<|j|<J$ with $\delta_0=1$, and $\delta_{-1}=-1$, 
by Lemma 3 there are $\sim p/2^{2J-2}$ values of $K$ with 
$\L Kp \L{K-j}p=\delta_j$ for each $j$.  Therefore 
(5.1) fails, for almost all $K$, by Proposition 3.
Appealing now to Theorem 2 we have proved that for 
$\sim \kappa_1p/2$  values of $K$ with $\L{K}{p}=-\L{K+1}{p}$,
the minimum of  $g_{p,K}(x)$  is  $<0$.  
For such $K$, $f_p(z)$ has two zeros on the arc from $\zeta_p^K$ 
to $\zeta_p^{K+1}$,  
so that the total number of such zeros is $\sim \kappa_1 p$.
Theorems 1 and 1$\frac{1}{2}$ follow.

\head 6.  Pseudo-Gauss Sums: Proof of the first part of Theorem 4\endhead

\noindent In this section, we wish to study the 
distribution of $f_p(\zeta^{K+\frac 12}_p)$.  By (2.3) and Proposition 4
we have (if $(\sqrt p >) J \to \infty$ as $p\to \infty$) for 
almost all $1\le K\le p-1$,
$$
\align
f_p(\zeta^{K+\frac 12}_p) &= \frac{if_p(\zeta_p)}{\pi}
\biggl(\sum_{|j| < J} \L{K-j}{p}\frac{1}{j+\frac 12} +
O\left(\frac{1}{J^{\frac 13}}\right)\biggr) \\
&= \eta \frac{\sqrt{p}}{\pi}
\biggl(\sum_{|j| < J} \L{K-j}{p}\frac{1}{j+\frac 12} +
O\left(\frac{1}{J^{\frac 13}}\right)\biggr),
\tag{6.1}\\
\endalign
$$
where $\eta=\pm 1$ or $\pm i$ is fixed.
Thus, by Lemma 3, we have that for any fixed real number  $\rho$
$$
\lim_{p\rightarrow \infty} \frac{1}{p} 
\#\biggl\{K : 1 \le K \le p \text{ and }H_p\left(\frac{K+\frac 12}{p}\right) 
< \rho \sqrt{p}\biggr\}
$$
exists and equals
$$
\lim_{J\rightarrow\infty} \ \text{Prob}\biggl(
\sum_{|j| < J} \frac{\delta_j}{j+\frac 12} < \pi \rho :\
\delta\in \Delta_{0,J} \biggr). \tag{6.2}
$$
(using the notation $\Delta_{J,K}$ of section 4).
One may obtain an expression for this probability as follows:  
Recall that $\int_0^\infty \frac{\sin y}{y} dy= \frac {\pi}{2}$, 
and so for any $k\neq 0$
$$
\frac{2}{\pi} \int_0^\infty \frac{\sin (kx)}{x} dx = \text{sgn}(k) 
\frac{2}{\pi} \int_0^{\infty} \frac{\sin(|k|x)}{x} dx 
= \text{sgn}(k) \frac{2}{\pi} \int_0^{\infty} \frac{\sin y}{y} dy
=\text{sgn}(k),
$$
where $\text{sgn}(k)$ is the sign of $k$ ($=1$ if $k>0$ and $-1$ if $k<0$).
Hence the probability (6.2) equals 
$$
\align
\frac{1}{2^{2J-1}} & \sum\Sb \delta\in \Delta_{0,J}\endSb
\biggl(\frac{1}{2} - \frac{1}{\pi} \int^{\infty}_{0}
\sin\biggl(\biggl(\sum_{|j|<J} \frac{\delta_j}{j+\frac 12} -
\pi \rho\biggr)x\biggr) \frac{dx}{x} \biggr) \\
&= \frac{1}{2} - \frac{1}{\pi} \int^{\infty}_{0} 
\frac{1}{2^{2J-1}}  \sum\Sb \delta\in \Delta_{0,J}\endSb
\biggl( \frac{e^{ix\left(\sum_{|j|<J} \frac{\delta_j}{j+1/2} - \pi \rho\right)}
- e^{-ix\left(\sum_{|j|<J} \frac{\delta_j}{j+1/2} -
\pi \rho\right)}} {2i} \biggr) \frac{dx}{x} \\
&= \frac{1}{2} - \frac{1}{\pi} \int^{\infty}_{0} \prod_{|j|<J}
\biggl(\frac{e^{\frac{ix}{j+1/2}} + e^{-\frac{ix}{j+1/2}}}{2}\biggr)
\left(\frac{e^{-ix\pi \rho} - e^{ix\pi \rho}}{2i}\right)\frac{dx}{x} \\
&= \frac{1}{2} + \frac{1}{\pi} \int^{\infty}_{x=0} \sin(\rho \pi x)
\prod_{|j|<J} \cos\left(\frac{2x}{2j+1}\right)\frac{dx}{x} .
\endalign
$$
Letting  $J \rightarrow \infty$, we get
$$
c_\rho = \frac{1}{2} + \frac{1}{\pi} \int^{\infty}_{0} \sin(\rho \pi x)
C(x) \frac{dx}{x}\qquad
\text{where} \qquad C(x):= \prod\Sb n \ge 1 \\ n \text{ odd} \endSb
\cos^2\left(\frac{2x}{n}\right) ,
$$
and thus Theorem 4 is proved.
Note that this integral does converge: For any $x>0$ we have
$$
C(x) \ll \frac{1}{2^{\frac{3x}{\pi}}}
$$
since this estimate is trivial for $x\leq 1$, and otherwise
we note that $|\cos(\frac{2x}{n})| < \frac 12$ if  
$3x/\pi < n < 6x/\pi$. Thus the part of the integral with $x\geq 1$ is
easily bounded.  Since $\sin(\rho\pi x)\ll \rho\pi x$, 
the portion of the integral from $0$ to $1$ is also easily bounded.  

\remark{Remark 1} 
We use the above to study the multiplicative average size of
$f_p(\zeta_p^{k+\frac 12})$.  Due to the symmetry of $c_\rho$ we have that 
$$
\frac{1}{p-1} \log \biggl( 
\prod_{k=1}^{p-1} \frac { f_p(\zeta_p^{k+\frac 12})} {\sqrt{p}} \biggr) 
= 2 \int_0^\infty \log \rho \ \text{d} (c_\rho -\tfrac 12) .
$$
Using our expression for $c_\rho$ one can show that this is
$$
=  \gamma + \log \pi -\
\int_0^1 \frac{C(x)-1}{x}\ dx \ -\ \int_1^\infty  \frac{C(x)}{x}\ dx .
$$
All of these integrals converge, though we do not know their exact values.
\endremark

\remark{Remark 2}
The expansion given in (6.1) for $f_p$, and the general technique involved,
is very similar to that used by Montgomery [5] in showing that, 
\roster
\item"{\ i)}" $|f_p(z)| \ll \sqrt{p} \log p$ for all $|z|=1$.
\item"{ii)}" If $p$ is sufficiently large then there exists some value
of $z$ with $|z|=1$ for which $|f_p(z)| > \frac{2}{\pi} \sqrt{p} \log\log p$.
\endroster
Indeed to prove a result like that in (ii) we note that we may select 
each $\delta_j$ equal to the sign of $j$ for $|j|<J=\varepsilon \log p$.
By Lemma 3 there are many such $K$ and we proceed as before with the
expansion in (6.1), but now taking a little more care over the
set of excluded $K$.
\endremark

\remark{Remark 3}
Fix $t\in (0,1)$. By the argument above, we have, for any fixed real number  $\rho$,
$$
\align
c_{\rho,t}:&=\lim_{p\rightarrow \infty} \frac{1}{p} 
\#\biggl\{K : 1 \le K \le p \text{ and }H_p\left(\frac{K+t}{p}\right) 
< \rho \sqrt{p}  \biggr\} \\
&= \lim_{J\rightarrow\infty} \ \text{Prob}\biggl(
\delta\in \Delta_{0,J}: \ \sum_{|j| < J} 
\frac{\delta_j}{j+t} < \frac{\pi \rho}{\sin(\pi t)} \biggr) \\
&=\frac{1}{2} + \frac{1}{\pi} \int^{\infty}_{x=0} \sin\left(\frac{\rho \pi x}{\sin(\pi t)}\right)
\prod_{j\in {\Bbb Z}} \cos\left(\frac{x}{j+t}\right)\frac{dx}{x} .\\
\endalign 
$$
\endremark

\remark{Remark 4} We can also use these techniques to investigate 
the distribution of values of $H_p(t)$ at $t=a/(p-1)$ for $1\leq a \leq p-1$. 
We note that if $K \sim \alpha p$ then 
$\zeta_{p-1}^K = \zeta_p^{K+\alpha} \{ 1 + o(\frac 1p) \}$. 
Therefore we can get
an expression similar to (6.1) for almost all $F_p(\zeta_{p-1}^K)$,
but now with $\sum_{|j|<J} \L{K-j}{p}  \frac{1}{j+\alpha}$ 
replacing the sum in (6.1), and multiplying the whole expression through 
by $\sin (\alpha \pi)$. Thus the density of those $K$, for which 
$H_p(K/(p-1)) \leq \rho \sqrt{p}$, is
$$
\frac{1}{2} + \frac{1}{\pi} \int^1_{\alpha=0} \int^{\infty}_{x=0} 
\sin\left( \frac{\rho \pi x}{\sin (\alpha \pi)} \right)
\prod\Sb m \in {\Bbb Z} \endSb
\cos \left(\frac{x}{m+\alpha}\right) \frac{dx}{x} d\alpha .
$$
We cannot see how to obtain a simpler expression.

It is not hard to modify this technique to determine the distribution
of values of the Fekete polynomial (or, in fact, $H_p(t)$) at any
``reasonably'' distributed set of values.

\endremark

\head 7. The distribution of $g(\tfrac 12)$ for $g \in 
\Cal F_J$ as $J\to \infty$. \endhead

\noindent We now look at the limiting distribution of $g(\frac 12)-4$ for
$g\in \Cal F_J$ as $J\to \infty$. 
Define, for $N\ge 1$,
$$
S_N(\d) = \sum_{|j+\frac 12| > N} \frac{\delta_j}{j+\frac 12},
$$
where each $\delta_j=1$ or $-1$ with probability $\frac 12$.
We will prove that the distribution function of $S_1(\d)$ decays
{\sl  double exponentially}.

\proclaim{Theorem 5} As $x\to \infty$, we have
Prob$(|S_1(\d)|>x) = \exp(-e^{\frac x2+O(1)})$.
\endproclaim

\demo{Proof of second part of Theorem 4} Note that 
Prob$(S_1(\d)>x) =$Prob$(S_1(\d)<-x) = \exp(-e^{\frac x2+O(1)})$, by symmetry.
Taking $x=\pi \rho$, the result follows from (6.2).
\enddemo

To prove Theorem 5 we study the $2k$-th moment of 
$S_N(\d)$, call it $M_N(k)$, 
that is, the expectation of $S_N(\d)^{2k}$. For example 
$$ 
M_N(1) = \sum_{|j+\frac 12|>N} \frac{1}{(j+\frac 12)^2}.
$$
Our aim is to determine 
the asymptotic behaviour of $M_1(k)$ for large $k$.

\proclaim{Proposition 5} For large $k$,
$$
M_1(k) = (2\log k - 2\log \log k + O(1))^{2k}.
$$
\endproclaim 

\demo{Proof} To establish the lower bound, consider $\d$ such 
that $\delta_j=1$ for all $1 \le |j+\frac 12| \le k/\log k$; and such that 
$S_{k/\log k}(\d) >0$.  The probability of this happening is  
$\asymp 1/2^{2k/\log k}$,  and
$S_1(\d) \ge 2\log k-2\log \log k+O(1)$  for such $\d$.  Hence
$$
M_1(k) \gg \frac{1}{2^{2k/\log k}} (2\log k-2\log \log k+O(1))^{2k}
= (2\log k - 2\log \log k + O(1))^{2k}.
$$

Now 
$$
M_N(k) = \sum_{j_1,j_2,\dots j_{2k}}
{\Bbb E}\left(  \frac{\delta_{j_1}}{j_1+\frac 12}
\frac{\delta_{j_2}}{j_2+\frac 12} \dots \frac{\delta_{j_{2k}}}{j_{2k}
+\frac 12} \right),
$$
where ${\Bbb E}$ stands for the expectation.  
Observe that a summand above is non-zero only 
if each value of $j$ appears an even
number of times amongst $j_1,j_2,\dots j_{2k}$. In particular $j_\ell=j_1$
for some $\ell>1$, and then 
${\Bbb E}(\prod_{1\leq i\leq 2k} \delta_{j_i}) ={\Bbb E}
(\prod_{1\leq i\leq 2k,\
i\ne 1,\ell} \delta_{j_i})$.  
Summing over all $2k-1$ possibilities for $\ell$
in the above, we deduce that 
$$
M_N(k) \le (2k-1) \sum_{|j+\frac 12|>N} \frac{1}{(j+\frac 12)^2} 
M_N(k-1), \tag{7.1}
$$
for all $k\ge 1$ and all $N\ge 1$. Iterating this inequality, we obtain
$$
\align
M_N(k) &\le (2k-1)\cdot(2k-3)\cdots 3\cdot 1 \cdot \biggl( \sum_{|j+\frac 
12|>N}
\frac{1}{(j+\frac 12)^2}\biggr)^k \\ &\le \frac{(2k)!}{k!2^k} 
\biggl(\frac{2}{N-\frac 12}\biggr)^k  = \frac{(2k)!}{k! (N-\frac 12)^k}.
\tag{7.2} \\
\endalign
$$

Now 
$$
|S_1(\d)-S_N(\d)| \le 2\lambda_N, \ \ \text{\rm where}\ 
\ 
\lambda_N:=\sum_{N\ge j+\frac 12 \ge 1} \frac{1}{j+\frac 12} 
= \log N+O(1) .
$$
Evidently the odd moments of $S_N(\d)$ are zero. 
Therefore, by the binomial theorem and (7.2), 
$$
\align
M_1(k) &= \sum_{j=0}^{k} \binom{2k}{2j}\  M_N(j) \ {\Bbb E}
(|S_1(\d)-S_N(\d)|^{2k-2j})\\
&\le \sum_{j=0}^{k} \binom{2k}{2j} \frac{(2j)!}{j!  (N-\frac 12)^{j}} 
\ (2\lambda_N)^{2k-2j} \\
&\le (2\lambda_N)^{2k} \sum_{j=0}^{k} \frac{1}{j!} \left( \frac{k^2}
{ (N-\frac12)\lambda_N^2} \right)^j 
\le (2\lambda_N)^{2k} \exp\biggl(\frac{k^2}{(N-\frac 12)\lambda_N^2}
\biggr).\\
\endalign
$$
Taking $N=k/\log k$ we obtain the upper bound of the Proposition.
\enddemo

\demo{Proof of Theorem 5} Take $k=c_1xe^{x/2}+O(1)$ for some $c_1>0$, and then

\noindent
Prob$(|S_1(\d)|>x)\leq  x^{-2k} M_1(k) \ll \exp(-c_2e^{x/2})$
for some constant $c_2>0$, if $c_1$ is sufficiently small, by Proposition 5.

The lower bound is more involved. Select integer $k$ so that 
$2\log k-2\log\log k$ is as close as possible to $x$. The 
contribution to $M_1(k)$ of those $\d$ with $|S_1(\d)|<x-c_3$ is
$\leq (x-c_3)^{2k}\leq M_1(k)/4$ if $c_3$ is sufficiently large. The 
contribution to $M_1(k)$ of those $\d$ with $|S_1(\d)|>x+c_3$ is
$\leq \int_{t>x+c_3} $Prob$(|S_1(\d)|>t) t^{2k} dt\ll \int_{t>x+c_3} 
\exp(-c_2e^{t/2}) t^{2k} dt \leq M_1(k)/4$ if $c_3$ is sufficiently large,
using the upper bound from the paragraph above. Thus
$M_1(k)/2 \leq $Prob$(x-c_3\leq |S_1(\d)|\leq x+c_3) (x+c_3)^k$
which implies that 
Prob$(|S_1(\d)|\geq x-c_3) \geq M_1(k)/2(x+c_3)^k \gg \exp(-c_4e^{x/2})$
for some constant $c_4>0$, by Proposition 5. 
Replacing $x-c_3$ by $x$ gives the lower bound
and thus our result.
\enddemo

\remark{Remark} We follow up on remark 3 of section 6.  The arguments
above (Theorem 5 and Proposition 5) hold just as well with ``$1/2$'' replaced by any fixed $t\in (0,1)$. 
Thus $1-c_{\rho,t}$ and $c_{-\rho,t}=\exp(-\exp( \pi \rho/2\sin(\pi t) +O(1)))$
for $\rho>0$.
\endremark

\head 8. Bounds on $\kappa_1$ \endhead

Applying the method of section 6, we note that for any real $\lambda$,
$$
\align
\pi_\lambda:&=\lim_{J\to \infty} \text{\rm Prob}\{g \in \Cal F_J : g(1/2)<4\lambda\} \\
&=\frac{1}{2} - \frac{1}{\pi} \int^{\infty}_0 \sin ( (1-\lambda) x)
\prod\Sb n \ge 3 \\ n \text{ odd} \endSb
\cos^2\left(\frac{x}{2n}\right) \frac{dx}{x} . \tag{8.1} \\
\endalign
$$
We can use this  to obtain  numerical bounds on  $\kappa_1$
using the following result.

\proclaim{Proposition 6}  We have
$\pi_{.013496\dots} \geq \kappa_1 \geq \pi_0$.
\endproclaim

Using Simpson's rule to compute the integrals in (8.1) we obtain
$.000813>\pi_{.013496\dots} \geq \kappa_1 \geq \pi_0>.000668$, from
which we deduce the bounds on $\kappa_0$ in the introduction.

\demo{Proof} Again selecting $x_0$ so that $g(x_0)$ is minimal, we have, by
definition, that
$$
\kappa_1 = \lim_{J\to \infty} \text{\rm Prob}\{g \in \Cal F_J : g(x_0)\leq 0\} .
$$
Since $g(x_0)\leq g(1/2)$ we deduce the lower bound on $\kappa_1$ above.

To get the upper bound, write  $x_0=\frac 12+\nu$ so that $|\nu |<\frac 12$.  
If $g(x_0)\le 0$ then
$$
\align
g(\tfrac 12)&\le g(\tfrac 12)-g(x_0) = 4-\frac{1}{x_0}-\frac{1}{1-x_0} + 
\sum\Sb |j|< J\\ j\neq 0,-1\endSb \frac{\delta_j (x_0-\frac 12)}
{(j+\frac 12)(j+x_0)}\\
&\le -\frac{4\nu^2}{\frac 14-\nu^2} + \sum_{j=1}^{\infty} 
\frac{|\nu|}{(j+\frac 12)(j+\frac 12+\nu)} 
+ \sum_{j=-\infty}^{-2} \frac{|\nu|}{(j+\frac 12)
(j+\frac 12+\nu)}\\
&=-\frac{4\nu^2}{\frac 14-\nu^2} + \sum_{j=1}^{\infty} \frac{2|\nu|}
{(j+\frac 12)^2-\nu^2}=-\frac{(2|\nu|+4\nu^2)}{\frac 14-\nu^2} 
+\pi \tan(\pi |\nu|).\\
\endalign
$$
Using Maple to compute the $\max_{\nu}$, we obtain
$$
g(\tfrac 12) \le \max_{|\nu|\le \frac 12} 
\biggl(\pi \tan(\pi |\nu|)-\frac{(2|\nu|+4\nu^2)}{\frac 14-\nu^2} 
\biggr) = 0.053986\dots ,
$$
the maximum being attained at $\nu = \pm .057052\dots$.
\enddemo

\remark{Remark}
One can refine the above to get better bounds for $\kappa_1$.  First note that $g(x)=1/x+1/(1-x)$ is the only element in $\Cal F_1$, and
in this case $x_0=1/2$; thus ``$1/2$'' appears in the definition of
$\pi_\lambda$. More generally, let $J$ be some positive integer. For each
$\gamma\in \Cal F_J$ select $\chi_0$ so that $\gamma(\chi_0)$ is minimal.
We again have $g(x_0)\leq g(\chi_0)$, so if $g(\chi_0)\leq 0$ then 
$g(x_0)\leq 0$. On the other hand, if $g(x_0)\le 0$ then
we can again get an explicit upper bound on $g(\chi_0)$ and proceed as above.
This can be used to give another proof that $\kappa_1$ exists.
\endremark

\head 9. Zeros off the unit circle \endhead

\demo{Proof of Theorem 3} 
Theorem 3 holds trivially if there is a zero of $f_p(t)$ on the unit circle in the arc from $\zeta_p^K$ to $\zeta_p^{K+1}$. Thus we shall henceforth assume 
that there is no such zero. Let $h(x):= H_p((K+x)/p) \big/ H_p(K/p)$, so that 
$|h(x)|=|f_p(\zeta_p^{K+x})/\sqrt{p}|$, and $h(x)$ is a continuous real-valued
function. Now the hypothesis implies that 
$h(y)<\epsilon$ for some $y\in (0,1)$ (in fact, $t=\zeta_p^{K+y}$), while our assumption above implies
that $h(x)\ne 0$ for all $x\in (0,1)$.  By (2.3) we have, uniformly for 
$|x|\leq 2/3$,
$$
\align
h(x)  &=\frac{\sin (\pi x)} {p} \left( \frac{1}{\sin (\pi x/p)} + 
\left( \frac Kp \right) \sum_{1\leq |K-k| <p/2} \frac{(k/p)}{\sin (\pi (x+K-k)/p)} \right) \\
&=1-(C+O(1))x, \quad \text{\rm where} \ C:= -\left( \frac Kp \right) \sum_{1\leq |K-k| <p/2} \frac{(k/p)}{K-k} .\\
\endalign
$$
So if $h(y)<\epsilon$ for some sufficiently small $y$ then $h(2y)=2h(y)-1+O(y)<0$,
contradicting our assumption. Therefore we may assume that $y\gg 1$, and
also $1-y\gg 1$ by the symmetric argument. Thus
$g_{p,K}(y)\ll \sqrt{p} |f_p(t)|/\sin(\pi y)\ll \epsilon p$ by (2.3), so that
$$g_{p,K}(x_0) \le g_{p,K}(y)\ll  \epsilon p$$ 
where $x_0$ is defined as in section 3.

Let $x_1 = x_0 - \epsilon^{1/2}$, and $x_2 = x_0 +\epsilon^{1/2}$, and then $\alpha_j=\zeta_p^{x_j}$ for $j=1,2$.  
Let $R=1-\epsilon^{\frac 13}/p$.  We shall
consider the variation in argument of 
$$
G(z):= i\L{K}{p} \frac{p}{f_p(\zeta_p)} \frac{f_p(z)}{z^p-1} = 
i \L{K}{p} \sum_{|K-k| < \frac p2} \L{k}{p} \frac{1}{z\zeta_p^{-k}-1},
$$
as $z$ goes around (in the anti-clockwise direction) 
the box bounded by the four curves, ${\Cal C}_1$, the arc of
the unit circle from $\alpha_1$ to $\alpha_2$, then ${\Cal C}_2$,
the straight line segment from  $\alpha_2$ to $R\alpha_2$, then 
${\Cal C}_3$, the arc of
the circle of radius $R$, from $R\alpha_2$ to $R\alpha_1$, then finally
${\Cal C}_4$,
the straight line segment from  $R\alpha_1$ back to $\alpha_1$.

We know that
$G(z)$ is real valued and positive on the arc ${\Cal C}_1$.  We shall show 
that $G(z)$ has positive imaginary part on ${\Cal C}_2$,
that $G(z)$ has negative real part on ${\Cal C}_3$, and
that $G(z)$ has negative imaginary part on ${\Cal C}_4$,
This shows that the change in argument of $G(z)$ is $2\pi$ as we go around our box, so that there is 
exactly one zero in our box. This implies a little more than Theorem 3.

To estimate $H(r,x):=G(r\zeta_p^{(K+x)/p})$ when $R\leq r\leq 1$, for
a value of $x\in [x_1,x_2]$,  we calculate the Taylor
series expansion around $r=1$, which is
$$
H(r,x) = g_{p,K}(x) - \frac{(1-r)^2}{2r}\biggl(\frac{p}{2\pi}
\biggr)^2 g_{p,K}^{\prime \prime}(x) + i \frac{1-r^2}{2r}\frac{p}{2\pi} g_{p,K}^{\prime}(x) + O\biggl(\frac{(1-r)^3}{r}p^4\biggr) .
$$
From the proof of Proposition 2 we have, since $x$ is bounded away from 0 and 1, $$
g_{p,K}(x) = g_{p,K}(x_0) + O((x-x_0)^2 p), \quad g_{p,K}^{\prime}(x) \asymp (x-x_0) p,
\quad \text{and} \ \ g_{p,K}^{\prime\prime}(x) \asymp p.
$$
Therefore 
$$
\align
\text{Im}(G(z)) = \text{Im}(H(r,x)) &\asymp \epsilon^{\frac{1}{2}} p^2 (1-r) + O ((1-r)\epsilon^{\frac{2}{3}} p^2) >0 \quad \text{\rm on} \ {\Cal C}_2, \\
\text{Im}(G(z)) = \text{Im}(H(r,x)) &\asymp - \epsilon^{\frac 12}p^2 (1-r) + O((1-r)\epsilon^{\frac 23} p^2) <0  \quad \text{\rm on} \ {\Cal C}_4, \\
\text{Re}(G(z)) = \text{Re}(H(r,x)) &\asymp - \epsilon^{\frac 23} p + O(\epsilon p)  <0 \quad \text{\rm on} \ {\Cal C}_3, \\
\endalign
$$
as required.
\enddemo
 
\remark{Remark} By (9.1) we see that 
$$
\max_{|z|=1}  \ |f_p(z)| \asymp \sqrt{p} \max_{K\in {\Bbb Z}}
\sum_{j\ne 0}  \frac{1}{j}  \left( \frac{K+j}{p} \right) .
$$
This again allows us to recover the results of Montgomery [5], as in remark 2
of section 6.
\endremark


\Refs
\ref\no 1
\by R.C. Baker and H.L. Montgomery
\paper Oscillations of Quadratic $L$-functions 
\inbook Analytic Number Theory (ed. B.C. Berndt et.al.)
\publ Birkh\"auser \publaddr Boston \yr 1990 \page 23--40 
\endref
\ref\no 2  
\by H. Davenport 
\book Multiplicative Number Theory {\rm (2nd ed.)}
\publ Springer-Verlag \publaddr New York \yr 1980
\endref
\ref\no 3
\by P. Erd\H os and P. Tur\'an
\paper On the distribution of roots of polynomials 
\jour Ann.~of Math.
\vol 51 \yr 1950 \pages 105--119
\endref
\ref\no 4
\by M. Fekete and G. P\' olya
\paper \" Uber ein Problem von Laguerre
\jour Rend.~Circ.~Mat.~Palermo \vol 34 \yr 1912 \page 89--120
\endref
\ref\no 5
\by  H.L. Montgomery
\paper An exponential polynomial formed with the Legendre symbol
\vol 37 \yr 1980 \jour Acta Arithmetica \page 375--380
\endref
\ref\no 6
\by G. P\'olya
\paper Verschiedene  Bemerkung zur Zahlentheorie
\jour Jber. deutsch Math. Verein \vol 28 \yr 1919 \page 31--40
\endref
\ref\no 7
\by M. Sambandham and V. Thangaraj
\paper On the average number of real zeros of a random 
trigonometric polynomial
\jour J. Indian Math. Soc \vol 47 \yr 1983 \page 139--150
\endref
\ref\no 8
\by A. Weil
\paper Sur les fonctions alg\' ebriques \` a corps de constantes fini
\jour C.R.~Acad.~Sci., Paris \vol 210 \yr 1940 \page 592--594
\endref
\endRefs
\enddocument